\documentstyle[11pt,amssymb,amsfonts]{article}

\begin{document}
\newcommand{\p}{\parallel }
\makeatletter \makeatother
\newtheorem{th}{Theorem}[section]
\newtheorem{lem}{Lemma}[section]
\newtheorem{de}{Definition}[section]
\newtheorem{rem}{Remark}[section]
\newtheorem{cor}{Corollary}[section]
\renewcommand{\theequation}{\thesection.\arabic {equation}}

\title{{\bf An equivariant Kastler-Kalau-Walze type
theorem} }
\author{  Yong Wang \\}

\date{}
\maketitle

\begin{abstract} In this paper, we prove an equivariant Kastler-Kalau-Walze type
theorem for spin manifolds without boundary. For $6$ dimensional spin manifolds with boundary, we also give an equivariant Kastler-Kalau-Walze type
theorem. Then we generalize this theorem to the general $n$ dimensional manifold. An equivariant Kastler-Kalau-Walze type theorem with torsion is also proved.
\\


\end{abstract}

\section{Introduction}
\quad The noncommutative residue found in [Gu] and [Wo] plays a
prominent role in noncommutative geometry. Several
years ago, Connes made a challenging observation that the
noncommutative residue of the square of the inverse of the Dirac
operator was proportional to the Einstein-Hilbert action, which we
call the Kastler-Kalau-Walze theorem. In [K], Kastler gave a
brute-force proof of this theorem. In [KW], Kalau and Walze proved
this theorem in the normal coordinates system simultaneously. In
[A], Ackermann gave a note on a new proof of this theorem
by means of the heat kernel expansion.\\
\indent On the other hand, Fedosov et al defined a noncommutative
residue on Boutet de Monvel's algebra and proved that it was a
unique continuous trace in [FGLS]. In [Wa1] and [Wa2], we
generalized some results in [C1] and [U] to the case of manifolds
with boundary . In [Wa3], We proved a Kastler-Kalau-Walze type
theorem for the Dirac operator and the signature operator for
$3,4$-dimensional manifolds with boundary. Recently, Ponge defined
lower dimensional volumes of Riemannian manifolds by the Wodzicki
residue in [Po]. In [Da], an equivariant residue was defined. So it is
nature question that whether we have an equivariant version of KKW theorem or not. \\
\indent  Atiyah-Singer index theorem is a very famous theorem. Its local version may be found in [BGV]. Considering a group action, Bismut defined the equivariant Bismut Laplacian instead of the Dirac Laplacian
and proved the infinitesimal equivariant index theorem. Motivated by Bismut's proof, in this paper we also use the equivariant Bismut Laplacian instead of the Dirac Laplacian to give an equivariant Kastler-Kalau-Walze type
theorem. Nextly, For $6$ dimensional spin manifolds with boundary, we also give an equivariant Kastler-Kalau-Walze type 
theorem. Then we generalize this theorem to the general $n$ dimensional manifold. An equivariant Kastler-Kalau-Walze type theorem with torsion is also proved. In Section 2, we give an equivariant Kastler-Kalau-Walze type
theorem. In Section 3, For $6$ dimensional spin manifolds with boundary, we also give an equivariant Kastler-Kalau-Walze type
theorem. In Section 4, we prove an equivariant Kastler-Kalau-Walze type
theorem for a $n$ dimensional spin manifold with boundary. In Section 5, an equivariant Kastler-Kalau-Walze type
theorem with torsion is given.

\section{An equivariant Kastler-Kalau-Walze type
theorem.}

   \quad Let $M$ be a smooth compact Riemannian $n$-dimensional
 manifold without boundary and $V$ be a vector bundle on $M$. Recall that a differential operator $P$
 is of Laplace type if it has locally the form
 $$P=-(g^{ij}\partial_i\partial_j+A^i\partial_i+B),\eqno(2.1)$$
 where $\partial_i$ is a natural local frame on $TM$ and $g_{i,j}=g(\partial_i,\partial_j)$ and
 $(g^{ij})_{1\leq i,j\leq m}$ is the inverse matrix associated
 to the metric matrix $(g_{i,j})_{1\leq i,j\leq m}$ on $M$, and $A^i$ and $B$ are smooth sections of
 ${\rm End}(V)$ on $M$ (endomorphism). If $P$ is a Laplace type operator
 of the form (2.1), then (see [Gi]) there is a unique connection
 $\nabla$ on $V$ and an unique endomorphism $E$ such that
 $$P=-[g^{ij}(\nabla_{\partial_i}\nabla_{\partial_j}-\nabla_{\nabla^L_{\partial_i}{\partial_j}})+E],\eqno(2.2)$$
 where $\nabla^L$ denotes the Levi-civita connection on $M$.
 Moreover (with local frames of $T^*M$ and $V$), $\nabla_{\partial_i}=
 \partial_i+\omega_i$ and $E$ are related to $g^{ij},~A^i$ and $B$
 through
 $$\omega_i=\frac{1}{2}g_{ij}(A^j+g^{kl}\Gamma_{kl}^j{\rm
 Id}),\eqno(2.3)$$
 $$E=B-g^{ij}(\partial_i(\omega_j)+\omega_i\omega_j-\omega_k\Gamma^k_{ij}),\eqno(2.4)$$
where $\Gamma_{ij}^k$ are the Christoffel coefficients of
$\nabla^L.$\\
   \indent Now we let $M$ be a $m$-dimensional oriented spin manifold
with Riemannian metric $g$. We recall that the Dirac operator $D$ is
locally given as follows in terms of orthonormal frames $e_i,~1\leq
i\leq n$ and natural frames $\partial_i$ of $TM$: one has
$$D=\sum_{i,j}g^{ij}c(\partial_i)\nabla^S_{\partial_j}=\sum_{i}c(e_i)\nabla^S_{e_i},\eqno(2.5)$$
where $c(e_i)$ denotes the Clifford action which satisfies the relation\\
 $$c(e_i)c(e_j)+c(e_j)c(e_i)=-2\delta_i^j,$$ and
$$\nabla^S_{\partial_i}=\partial_i+\sigma_i,~~\sigma_i=\frac{1}{4}\sum_{j,k}\left<\nabla^L_{\partial_i}e_j,e_k\right>c(e_j)c(e_k).\eqno(2.6)$$
Let
$$\partial^j=g^{ij}\partial_i,~~\sigma^i=g^{ij}\sigma_j,~~\Gamma^k=g^{ij}\Gamma_{ij}^k.\eqno(2.7)$$
By (6a) in [Ka], we have
$$D^2=-g^{ij}\partial_i\partial_j-2\sigma^j\partial_j+\Gamma^k\partial_k-g^{ij}[\partial_i(\sigma_j)+\sigma_i\sigma_j-\Gamma_{ij}^k\sigma_k]+\frac{1}{4}r,
\eqno(2.8)$$ where $r$ is the scalar curvature.\\
\indent Let a compact group $G$ act isometrically on $M$ and preserve the spin structure. This action generates a Killing vector field $X$. Let $L_X$ be the Lie derivation on the Spinors bundle.
The Levi-Civita connection $\nabla^L$ lifts a Clifford connection $\nabla^S$. Following Bismut, we define the equivariant Bismut Laplacian.
$$H_X=(D+\frac{1}{4}c(X))^2+L_X; L_X=\nabla^S_X+\mu(X),\eqno(2.9)$$
So$$H_X=D^2+\frac{1}{4}Dc(X)+\frac{1}{4}c(X)D+\nabla^S_X+\mu(X)-\frac{1}{16}|X|^2.\eqno(2.10)$$
Let $X=X_j\partial_j$. So $L_X=X_j\partial_j+\frac{1}{4}X_j\sigma_j+\mu(X).$ Then by (10) in [Wa5], we have

$$H_X=-g^{ij}\partial_i\partial_j+[X_j-2\sigma^j+\Gamma^j+\frac{1}{4}c(\partial^j)c(X)+\frac{1}{4}c(X)c(\partial^j)]\partial_j$$
$$+g^{ij}[-\partial_i(\sigma_j)-\sigma_i\sigma_j+\Gamma_{ij}^k\sigma_k+\frac{1}{4}c(\partial_i)\partial_j(c(X))+\frac{1}{4}c(\partial_i)\sigma_jc(X)$$
$$+\frac{1}{4}c(X)c(\partial_i)\sigma_j]+\frac{1}{4}r-\frac{1}{16}|X|^2+\frac{1}{4}X_j\sigma_j+\mu(X).\eqno(2.11)$$
$$(D+\frac{1}{4}c(X))^2=-[g^{ij}(\nabla_{\partial_i}\nabla_{\partial_j}-\nabla_{\nabla^L_{\partial_i}{\partial_j}})+E],\eqno(2.12)$$
$$(D+\frac{1}{4}c(X))^2+L_X=-[g^{ij}(\widetilde{\nabla}_{\partial_i}\widetilde{\nabla}_{\partial_j}-\widetilde{\nabla}_{\nabla^L_{\partial_i}{\partial_j}})+\widetilde{E}],\eqno(2.13)$$
By (2.2)-(2.4), we have
$$\widetilde{\omega_i}=\omega_i+\frac{1}{2}g_{ij}X_j, \widetilde{A}_j=A_j+X_j,\widetilde{B}=B+\frac{1}{4}X_j\sigma_j+\mu(X).\eqno(2.14)$$
By (11) in [wa5] and (2.2)-(2.4), we have
$$\widetilde{E}=E+\frac{1}{4}X_j\sigma_j+\mu(X)-g^{ij}[\partial(\frac{1}{2}g_{jk}X_k)+\frac{1}{2}g_{il}X_l\omega_j+\frac{1}{4}g_{ik}g_{jl}X_kX_l-\frac{1}{2}g_{kl}X_l\Gamma^k_{ij}].\eqno(2.15)$$
$$E=-\frac{1}{4}r+\frac{1}{16}|X|^2+\frac{1}{2}[e_j(\frac{1}{4}c(X))c(e_j)-c(e_j)e_j(\frac{1}{4}c(X))]+\sum_i<e_i,X>^2.\eqno(2.16)$$
By (2.14), we have
$${\widetilde{\nabla}}_Y={\nabla}_Y+\frac{3}{4}g(X,Y).\eqno(2.17)$$
By ${\rm tr}\sigma_j=0, {\rm tr}\sigma^j=0$ and direct computations, we get
$${\rm tr}(\frac{r}{6}+\widetilde{E})={\rm dim}(S(TM))\left\{-\frac{1}{12}r+\frac{1}{16}|X|^2+\sum_i<e_i,X>^2\right.$$
$$-g^{ij}\left[\partial_i(\frac{1}{2}g_{jk}X_k)+\frac{1}{2}g_{il}X_l[\frac{1}{2}g_{j\alpha}(\Gamma^\alpha+g^{kl}\Gamma^\alpha_{kl})]\right]$$

$$\left.-\frac{1}{4}\left<\partial^j,X\right>g_{ij}
+\frac{1}{4}g_{ik}g_{jl}X_kX_l-\frac{1}{2}g_{kl}X_l\Gamma^k_{ij}\right\}+{\rm tr}(\mu(X)).\eqno(2.18)$$
Let ${\rm dim}M=m$. By [A], we know that
$${\rm Wres}(H_X^{-\frac{m}{2}+1})=\frac{m-2}{(4\pi)^{\frac{m}{2}}\Gamma(\frac{m}{2})}{\rm tr}(\frac{r}{6}+\widetilde{E}),\eqno(2.19)$$
where Wres denotes the noncommutative residue (see [A]). So by (2.18) and (2.19), we have

\noindent{\bf Theorem 1} {\it The following equality holds}
$${\rm Wres}(H_X^{-\frac{m}{2}+1})=\frac{m-2}{(4\pi)^{\frac{m}{2}}\Gamma(\frac{m}{2})}\int_M
\left[{\rm dim}(S(TM))\left\{-\frac{1}{12}r+\frac{1}{16}|X|^2+\sum_i<e_i,X>^2\right.\right.$$
$$-g^{ij}\left[\partial_i(\frac{1}{2}g_{jk}X_k)+\frac{1}{2}g_{il}X_l[\frac{1}{2}g_{j\alpha}(\Gamma^\alpha+g^{kl}\Gamma^\alpha_{kl})]\right]$$
$$\left.\left.-\frac{1}{4}\left<\partial^j,X\right>g_{ij}
+\frac{1}{4}g_{ik}g_{jl}X_kX_l-\frac{1}{2}g_{kl}X_l\Gamma^k_{ij}\right\}+{\rm tr}(\mu(X))\right].\eqno(2.20)$$

\section{An equivariant KKW theorem for $6$-dimensional manifolds with boundary}

\quad Let $M$ be a $6$-dimensional compact oriented spin manifold with boundary $\partial M$. We assume that the metric $g^M$ on $M$
has the following form near the boundary,
 \begin{equation}
 g^M=\frac{1}{h(x_n)}g^{\partial M}+dx_n^2,
\end{equation}
where $g^{\partial M}$ is the metric on  ${\partial M}$.
Let a compact group $G$ act isometrically on $M$ and preserve the spin structure. Near the boundary, the group action is not necessarily product action. Let $\widetilde{{\rm Wres}}$
denote the noncommutative residue for manifolds with boundary (see [FGLS]). In the following, we want to compute $\widetilde{{\rm Wres}}[\pi^+H_X^{-1}\circ\pi^+H_X^{-1}].$
By the definition of the noncommutative residue for manifolds with boundary, we have
$$\widetilde{{\rm Wres}}[\pi^+H_X^{-1}\circ\pi^+H_X^{-1}]=\int_M\int_{|\xi|=1}{\rm
trace}_{S(TM)}[\sigma_{-6}(H^{-2}_X)]\sigma(\xi)dx+\int_{\partial
M}\Phi,\eqno(3.1)$$ \noindent where
$$\Phi=\int_{|\xi'|=1}\int^{+\infty}_{-\infty}\sum^{\infty}_{j, k=0}
\sum\frac{(-i)^{|\alpha|+j+k+1}}{\alpha!(j+k+1)!}$$
$$\times {\rm trace}_{S(TM)}
[\partial^j_{x_n}\partial^\alpha_{\xi'}\partial^k_{\xi_n}
\sigma^+_{r}(H^{-1}_X)(x',0,\xi',\xi_n)\times
\partial^\alpha_{x'}\partial^{j+1}_{\xi_n}\partial^k_{x_n}\sigma_{l}
(H^{-1}_X)(x',0,\xi',\xi_n)]d\xi_n\sigma(\xi')dx',\eqno(3.2)$$
\noindent where the sum is taken over $
r-k-|\alpha|+l-j-1=-6,~~r\leq -2,l\leq -2$.
Interior term comes from Theorem 1. We only compute the boundary term.
By (2.11), we have the symbol $\sigma_{-2}(H_X)$ and  $\sigma_{-3}(H_X)$ of $H_X$
By the symbol composition formula, we have
$$1=\sum_{|\alpha|=0}^2\frac{1}{\alpha!}\partial_\xi^\alpha(\sigma(H_X))D_x^\alpha(\sigma(H_X^{-1}),\eqno(3.3)$$
where $D_x^\alpha=(-i)^{|\alpha|}\partial^\alpha_x$.
Let $\sigma(H_X)=p_2+p_1+p_0, \sigma(H^{-1}_X)=r_{-2}+r_{-3}+r_{-4}+\cdots,$ then we can get
$$r_{-2}(H^{-1}_X)=|\xi|^{-2}; p_1r_{-2}+p_2r_{-3}+\sum_j\partial_{\xi_j}D_{x_j}r_{-2}=0.\eqno(3.4)$$
By (2.11) and (3.4) and Lemma 1 in [Wa4], we get
$$r_{-3}(H^{-1}_X)=r_{-3}(D^{-2})-\sqrt{-1}|\xi|^4(X_j-\frac{1}{2}<X,\partial_j>)\xi_j,\eqno(3.5)$$
where $r_{-3}(D^{-2})$ in Lemma 1 in [Wa4]. \\
\quad
\indent Now we can compute $\Phi$, since the sum is taken over $
-r-l+k+j+|\alpha|=-5,~~r,l\leq-2,$ then we have the same five cases as in [Wa4].
By $r_{-2}(D^{-2})=r_{-2}(H^{-1}_X)$, we know that our $\overline{\rm cases}$ 1 (i) (ii) (iii) equal cases (i) (ii) (iii) in [Wa4]. So the sum of these three cases is zero by [Wa4].\\
\noindent  ${\bf \overline{case}}$~II)~$r=-2,~l=-3~k=j=|\alpha|=0.$\\
\noindent By (3.2) and an integration by parts and (19) in [Wa4], we get
$${ \overline{\rm case}
~II)}=
-\sqrt{-1}\int_{|\xi'|=1}\int^{+\infty}_{-\infty}
{\rm trace} [\pi^+_{\xi_n}\sigma_{-2}(H^{-1}_X)\times
\partial_{\xi_n}\sigma_{-3}(H^{-1}_X)](x_0)d\xi_n\sigma(\xi')dx'$$
$$=\sqrt{-1}\int_{|\xi'|=1}\int^{+\infty}_{-\infty}
{\rm trace}\{\partial_{\xi_n}\pi^+_{\xi_n}\sigma_{-2}(H^{-1}_X)\times
[\sigma_{-3}(D^{-2})$$
$$-\sqrt{-1}|\xi|^4(X_j-\frac{1}{2}<X,\partial_j>)\xi_j]\}(x_0)d\xi_n\sigma(\xi')dx'.)$$
$$={\rm case~ II}+\sqrt{-1}\int_{|\xi'|=1}\int^{+\infty}_{-\infty}{\rm tr}\frac{i}{2(\xi_n-i)^2}(-\sqrt{-1})|\xi|^4(X_j-\frac{1}{2}<X,\partial_j>)\xi_j](x_0)d\xi_n\sigma(\xi')dx'\eqno(3.6)$$
where case II is in [Wa4].
By $\int_{|\xi'|=1}\xi_j\sigma(\xi')=0$ for $j<m$ and the metric on $M$ and the Cauchy integral formula, we get
$$\sqrt{-1}\int_{|\xi'|=1}\int^{+\infty}_{-\infty}\frac{i}{2(\xi_n-i)^2}(-\sqrt{-1})|\xi|^4(X_j-\frac{1}{2}<X,\partial_j>)\xi_j](x_0)d\xi_n\sigma(\xi')dx'$$
$$=-\frac{1}{32}\Omega_4X_m{\rm Vol}_{\partial M}\eqno(3.7)$$
where $\Omega_4$ is the canonical volume of the sphere $S^4$.\\
\noindent ${\bf  \overline{case}~ III)}$~$r=-3,~l=-2,~k=j=|\alpha|=0$\\
$${ \overline{\rm case}~ III)}=-i\int_{|\xi'|=1}\int^{+\infty}_{-\infty}
{\rm trace} [\pi^+_{\xi_n}\sigma_{-3}(H^{-1}_X)\times
\partial_{\xi_n}\sigma_{-2}(H^{-1}_X)](x_0)d\xi_n\sigma(\xi')dx'.\eqno(3.8)$$
Let $A=\sqrt{-1}|\xi|^4(X_j-\frac{1}{2}<X,\partial_j>)\xi_j$. So
$${ \overline{\rm case}~ III)}={{\rm case}~ III)}+\sqrt{-1}
\int_{|\xi'|=1}\int^{+\infty}_{-\infty}
{\rm trace}  [\pi^+_{\xi_n}A\partial_{\xi_n}\sigma_{-2}(H^{-1}_X)](x_0)d\xi_n\sigma(\xi')dx'.\eqno(3.9)$$
Similarly to (23) in [Wa4], using the Cauchy integral formula, we get
$${ \overline{\rm case}~ III)}={ {\rm case}~ III)}+\frac{1}{32}\Omega_4X_m{\rm Vol}_{\partial M}\eqno(3.10)$$
By [Wa4], we know that ${ {\rm case}~ II)}+{ {\rm case}~ III)}=0$, so by (3.6) (3.7) and (3.10), we get
$\Phi=0$
So we get the following theorem\\

\noindent{\bf Theorem 2} {\it For $6$ dimensional spin manifolds with boundary, the following equality holds}
$$\widetilde{{\rm Wres}}[\pi^+H_X^{-1}\circ\pi^+H_X^{-1}]=\frac{1}{4\pi^3}\int_M
\left[\left\{-\frac{1}{12}r+\frac{1}{16}|X|^2+\sum_i<e_i,X>^2\right.\right.$$
$$-g^{ij}\left[\partial_i(\frac{1}{2}g_{jk}X_k)+\frac{1}{2}g_{il}X_l[\frac{1}{2}g_{j\alpha}(\Gamma^\alpha+g^{kl}\Gamma^\alpha_{kl})]\right]$$
$$\left.\left.-\frac{1}{4}\left<\partial^j,X\right>g_{ij}
+\frac{1}{4}g_{ik}g_{jl}X_kX_l-\frac{1}{2}g_{kl}X_l\Gamma^k_{ij}\right\}+{\rm tr}(\mu(X))\right].\eqno(3.12)$$\\

\noindent{\bf Remark.} We may extend this theorem to the higher dimensional case.\\

In order to get the nonzero boundary term, we use a function $f$ on $M$ to perturb $\widetilde{{\rm Wres}}[\pi^+H_X^{-1}\circ\pi^+H_X^{-1}]$.
We find only the term $\overline{\rm case }$ I ii) change and other terms does not change. We get \\

\noindent{\bf Theorem 3} {\it For $6$ dimensional spin manifolds with boundary, the following equality holds}
$$\widetilde{{\rm Wres}}[\pi^+fH_X^{-1}\circ\pi^+H_X^{-1}]=\frac{1}{4\pi^3}\int_Mf
\left[\left\{-\frac{1}{12}r+\frac{1}{16}|X|^2+\sum_i<e_i,X>^2\right.\right.$$
$$-g^{ij}\left[\partial_i(\frac{1}{2}g_{jk}X_k)+\frac{1}{2}g_{il}X_l[\frac{1}{2}g_{j\alpha}(\Gamma^\alpha+g^{kl}\Gamma^\alpha_{kl})]\right]$$
$$\left.\left.-\frac{1}{4}\left<\partial^j,X\right>g_{ij}
+\frac{1}{4}g_{ik}g_{jl}X_kX_l-\frac{1}{2}g_{kl}X_l\Gamma^k_{ij}\right\}+{\rm tr}(\mu(X))\right]-\pi\Omega_4\int_{\partial_M}\partial_{x_n}f|_{x_n=0}d{\rm Vol}_{\partial_M}.\eqno(3.13)$$\\

\section{An equivariant general KKW theorem for manifolds with boundary}

\quad Let $M$ be a $n=\overline{n}+2$-dimensional compact oriented spin manifold with boundary $\partial M$ and $\overline{n}$ is an even integer. We will compute
$\widetilde{{\rm Wres}}[\pi^+H_X^{-1}\circ\pi^+H_X^{-\frac{\overline{n}}{2}+1}]$. By the definition of $\widetilde{{\rm Wres}}$ (see [FGLS]), we only compute the term
$$\Phi'=\int_{|\xi'|=1}\int^{+\infty}_{-\infty}\sum^{\infty}_{j, k=0}
\sum\frac{(-i)^{|\alpha|+j+k+1}}{\alpha!(j+k+1)!}$$
$$\times {\rm trace}_{S(TM)}
[\partial^j_{x_n}\partial^\alpha_{\xi'}\partial^k_{\xi_n}
\sigma^+_{r}(H^{-1}_X)(x',0,\xi',\xi_n)$$
$$\times
\partial^\alpha_{x'}\partial^{j+1}_{\xi_n}\partial^k_{x_n}\sigma_{l}
(H^{-\frac{\overline{n}}{2}+1})(x',0,\xi',\xi_n)]d\xi_n\sigma(\xi')dx',\eqno(4.1)$$
\noindent where the sum is taken over $
r-k-|\alpha|+l-j-1=-n,~~r\leq -2,l\leq 2-\overline{n}$.
Similar to Section 3 and [WW], we divide $\Phi'$ into five terms and the first three terms have the same expressions with the three term in [WW].\\
$\overline{\rm case} II: r=-2, l=1-\overline{n}, k=j=|\alpha|=0$\\
$${ \overline{\rm case}
~II)}=
-\sqrt{-1}\int_{|\xi'|=1}\int^{+\infty}_{-\infty}
{\rm trace} [\partial_{\xi_n}\pi^+_{\xi_n}\sigma_{-2}(H^{-1}_X)\times
\sigma_{1-\overline{n}}(H^{-\frac{\overline{n}}{2}+1})](x_0)d\xi_n\sigma(\xi')dx'\eqno(4,2)$$
Similar to (3.30) in [WW], we have
$$\sigma_{1-\overline{n}}(H^{-\frac{\overline{n}}{2}+1}_X)=\frac{\overline{n}-2}{2}\sigma^{-\frac{\overline{n}}{2}+2}_2\sigma_{-3}(H^{-1}_X)-\sqrt{-1}\sum_{k=0}^{\frac{\overline{n}}{2}-3}
\partial_{\xi_\mu}\sigma^{-\frac{\overline{n}}{2}+k+2}\partial_{x_\mu}\sigma^{-1}_2(\sigma^{-1}_2)^k.\eqno(4.3)$$
By (3.30) in [WW] and (3.5), we have
$$\sigma_{1-\overline{n}}(H^{-\frac{\overline{n}}{2}+1}_X)=\sigma_{1-\overline{n}}(D^{-\overline{n}+2})
-\sqrt{-1}|\xi|^{-4}(X_j-\frac{1}{2}<X,\partial_j>)\xi_j\frac{\overline{n}-2}{2}\sigma^{-\frac{\overline{n}}{2}+2},\eqno(4.4)$$
So
$$\overline{{\rm case}} II={\rm case} II+\sqrt{-1}\int_{|\xi'|=1}\int^{+\infty}_{-\infty}{\rm trace}\left[\frac{1}{2(\xi_n-i)^2}\right.$$
$$\left.\times|\xi|^{-4}(X_j-\frac{1}{2}<X,\partial_j>)\xi_j\frac{\overline{n}-2}{2}|\xi|^{-\overline{n}+4}\right]\sigma(\xi')dx':={\rm case} II+A,\eqno(4.5)$$
By $\int_{|\xi'|=1}\xi_j\sigma(\xi')=0$ for $j<m$ and the metric on $M$ and the Cauchy integral formula, we get
$$A=\frac{(2-\overline{n})2^{\frac{\overline{n}}{2}-2}}{(\frac{\overline{n}}{2}+1)!}
X_n\Omega(S_{\overline{n}})
\left[\frac{\xi_n}{(\xi_n+i)^{\frac{\overline{n}}{2}}}\right]^{(\frac{\overline{n}}{2}+1)}\left|_{\xi_n=i}\right..\eqno(4.6)$$\\

$\overline{case} III, r=-3,l=2-\overline{n}, k=j=|\alpha|=0$\\
$${ \overline{\rm case}~ III)}=-i\int_{|\xi'|=1}\int^{+\infty}_{-\infty}
{\rm trace} [\pi^+_{\xi_n}\sigma_{-3}(H^{-1}_X)\times
\partial_{\xi_n}\sigma_{2-\overline{n}}(H_X^{\frac{2-\overline{n}}{2}}))](x_0)d\xi_n\sigma(\xi')dx'.\eqno(4.7)$$
By (3.5) and (3.33) in [WW], we have
$$ \overline{{\rm case}}III={\rm case}III-\sqrt{-1}
\int_{|\xi'|=1}\int^{+\infty}_{-\infty}
{\rm trace} [\pi^+_{\xi_n}(-\sqrt{-1}|\xi|^{-4}(X_j-\frac{1}{2}<X,\partial_j>)\xi_j)$$
$$\times\partial_{\xi_n}\sigma_{2-\overline{n}}(H^{\frac{2-\overline{n}}{2}}_X)](x_0)d\xi_n\sigma(\xi')dx':={\rm case}III+B\eqno(4.8)$$
Similar to case II, through some computations, we get

$$B=\frac{\frac{\overline{n}}{2}-1}{(\frac{\overline{n}}{2}+1)!}
2^{\frac{\overline{n}}{2}+1}
X_n\Omega(S_{\overline{n}})
\left[\frac{\xi_n}{(\xi_n+i)^{\frac{\overline{n}}{2}}}\right]^{(\frac{\overline{n}}{2}+1)}\left|_{\xi_n=i}\right..\eqno(4.9)$$\\
By (4.5), (4.6), (4.8) and (4.9), we have
$$\Phi'=\Phi+A+B=\Phi+\frac{3n-6}{(\frac{\overline{n}}{2}+1)!}
2^{\frac{\overline{n}}{2}-2}
X_n\Omega(S_{\overline{n}})
\left[\frac{\xi_n}{(\xi_n+i)^{\frac{\overline{n}}{2}}}\right]^{(\frac{\overline{n}}{2}+1)}\left|_{\xi_n=i}\right..\eqno(4.10)$$\\
By Theorem 1 and (4.10) and (3.42) in [WW], we get\\

\noindent {\bf Theorem 4} {\it The following equality holds}
$$\widetilde{{\rm Wres}}[\pi^+H_X^{-1}\circ\pi^+H_X^{-\frac{\overline{n}}{2}+1}]=
\frac{n-2}{(4\pi)^{\frac{n}{2}}\Gamma(\frac{n}{2})}2^{\frac{n}{2}}\int_M
\left[\left\{-\frac{1}{12}r+\frac{1}{16}|X|^2+\sum_i<e_i,X>^2\right.\right.$$
$$-g^{ij}\left[\partial_i(\frac{1}{2}g_{jk}X_k)+\frac{1}{2}g_{il}X_l[\frac{1}{2}g_{j\alpha}(\Gamma^\alpha+g^{kl}\Gamma^\alpha_{kl})]\right]$$
$$\left.\left.-\frac{1}{4}\left<\partial^j,X\right>g_{ij}
+\frac{1}{4}g_{ik}g_{jl}X_kX_l-\frac{1}{2}g_{kl}X_l\Gamma^k_{ij}\right\}+{\rm tr}(\mu(X))\right]$$
$$-\frac{\overline{n}-2}{\overline{n}+1}\frac{\pi i}{(\frac{\overline{n}}{2}+2)!}
2^{\frac{\overline{n}}{2}-1}A_0\Omega(S_{\overline{n}})\int_{\partial_M}K{\rm dvol}_{\partial M}$$
$$+\frac{3n-6}{(\frac{\overline{n}}{2}+1)!}
2^{\frac{\overline{n}}{2}-2}
\Omega(S_{\overline{n}})
\left[\frac{\xi_n}{(\xi_n+i)^{\frac{\overline{n}}{2}}}\right]^{(\frac{\overline{n}}{2}+1)}\left|_{\xi_n=i}\right.\int_{\partial_M}X_n{\rm dvol}_{\partial M},\eqno(4.11)$$\\
{\it where $K$ is an extrinsic curvature and $A_0$ is a constant (see (3.42) in [WW]).}\\

\section {An equivariant KKW theorem with torsion}

 \quad Let $T$ be a real three form on $M$. Let $D_T=D+T$ where $T$ denotes the three form induces the Clifford action. Then $D_T$ is self adjoint operator and
 $$D^2_T=\triangle+\frac{1}{4}r+\frac{3}{2}dT-\frac{3}{4}||T||^2,\eqno(5.1)$$ 
 Define the Bismut Laplacian with torsion $H^T_X=(D_T+\frac{1}{4}c(X))^2+L_X$.
 So
 $$H^T_X=H_X+\frac{3}{2}dT-\frac{3}{4}||T||^2+\frac{1}{4}Tc(X)+\frac{1}{4}c(X)T,\eqno(5.2)$$
so by the formulas (2.3) and (2.4), we get
$$E_{H^T_X}=E_{H_X}+\frac{3}{2}dT-\frac{3}{4}||T||^2+\frac{1}{4}Tc(X)+\frac{1}{4}c(X)T,\eqno(5.3)$$
We get\\

\noindent{\bf Theorem 5} {\it The following equality holds}
$${\rm Wres}((H^T_X)^{-\frac{m}{2}+1})=\frac{m-2}{(4\pi)^{\frac{m}{2}}\Gamma(\frac{m}{2})}2^{\frac{m}{2}}\int_M
\left[\left\{-\frac{1}{12}r+\frac{1}{16}|X|^2+\sum_i<e_i,X>^2\right.\right.$$
$$-g^{ij}\left[\partial_i(\frac{1}{2}g_{jk}X_k)+\frac{1}{2}g_{il}X_l[\frac{1}{2}g_{j\alpha}(\Gamma^\alpha+g^{kl}\Gamma^\alpha_{kl})]\right]$$
$$\left.\left.-\frac{1}{4}\left<\partial^j,X\right>g_{ij}
+\frac{1}{4}g_{ik}g_{jl}X_kX_l-\frac{1}{2}g_{kl}X_l\Gamma^k_{ij}\right\}+{\rm tr}(\mu(X)
+\frac{3}{2}dT-\frac{3}{4}||T||^2+\frac{1}{2}Tc(X)
)\right].\eqno(5.4)$$\\

Bismut proved the local infinitesimal equivariant index formula by the Bismut Laplacian. Bismut also proved the local index theorem with torsion (see BGV, section 8.3 and [Bi]).
By above the the Bismut Laplacian with torsion, we may prove a local infinitesimal equivariant index formula with torsion as following (details will appear elsewhere).\\

\noindent{\bf Theorem 6} {\it Let $dT=0$ and $i_XT=0$ and $X$ is small, we have}
$${\rm lim}_{t\rightarrow 0}{\rm str}[{\rm exp}(-tH^T_X)(x,x)]d{\rm vol}_M=(2\pi\sqrt{-1})^{-\frac{n}{2}}\widehat{A}(F^{-T}_g(X)),\eqno(5.5)$$ 
{where $F^{-T}_g(X)=R^{-T}+\mu(X)$ and $\mu(X)$ is the moment map.}\\

\noindent{\bf Acknowledgement:}~~The work of the author was supported by NSFC. 11271062 and NCET-13-0721.\\ \\

\noindent{\bf References}\\

\noindent [A] T. Ackermann, {\it A note on the Wodzicki residue,} J.
Geom. Phys., 20, 404-406, 1996.\\
\noindent [BGV] N. Berline, E. Getzler, M. Vergne, {\it Heat Kernals
and Dirac Operators,}
Springer-Verlag, Berlin, 1992.\\
\noindent [Bi] J. Bismut, {\it A local index theorem for non-Kahler manifolds}. Math. Ann. 284 (1989), no. 4, 681-699.\\
\noindent [C1] A. Connes, {\it Quantized calculus and applications,}
XIth International Congress of Mathematical Physics (paris,1994),
15-36, Internat Press, Cambridge, MA, 1995.\\
\noindent[Da] S. Dave,{\it An equivariant noncommutative residue}. J. Noncommut. Geom. 7 (2013), no. 3, 709-735.\\
\noindent [FGLS] B. V. Fedosov, F. Golse, E. Leichtnam, and E.
Schrohe. {\it The noncommutative residue for manifolds with
boundary,} J. Funct.
Anal, 142:1-31,1996.\\
\noindent[Gi]P. B. Gilkey.: Invariance theory, the Heat equation, and the Atiyah-Singer Index theorem. Inc., USA, (1984).\\
 \noindent [Gu] V.W. Guillemin, {\it A new proof of Weyl's
formula on the asymptotic distribution of eigenvalues}, Adv. Math.
55 no.2, 131-160, 1985.\\
\noindent [K] D. Kastler, {\it The Dirac operator and gravitiation,}
Commun. Math. Phys, 166:633-643, 1995.\\
\noindent [KW] W. Kalau and M.Walze, {\it Gravity, non-commutative
geometry, and the Wodzicki residue,} J. Geom. Phys., 16:327-344, 1995.\\
\noindent[Po] R. Ponge, {\it Noncommutative geometry and lower
dimensional volumes in Riemannian geometry,} Lett. Math. Phys. 83
(2008), no. 1, 19--32.\\
\noindent [U] W. J. Ugalde, {\it Differential forms and the Wodzicki
residue,} J. Geom. Phys. 58 (2008), no. 12, 1739--1751.\\
\noindent [Wa1] Y. Wang, {\it Differential forms and the Wodzicki
residue for manifolds with boundary,} J. Geom. Phys.,
56:731-753, 2006.\\
 \noindent
[Wa2] Y. Wang, {\it Differential forms and the noncommutative
residue for manifolds with boundary in the non-product Case,} Lett.
math. Phys., 77:41-51, 2006.\\
\noindent [Wa3] Y. Wang, {\it Gravity and the noncommutative residue
for manifolds with boundary,} Lett. Math. Phys. 80 (2007), no. 1,
37--56.\\
\noindent [Wa4], Y. Wang, {\it Lower dimensional volume and Kastler-Kalau-Walze theorem for manifolds with boundary}, Commun. Theor. Phys., 54 (2010), 38-42.\\
\noindent [Wa5] Y. Wang, {\it A Kastler-Kalau-Walze type theorem and the spectral action for perturbations of Dirac operators on manifolds with boundary}. Abstr. Appl. Anal. 2014, Art. ID 619120, 13 pp.\\
\noindent [WW] J. Wang and Y. Wang, {\it A general Kastler-Kalau-Walze type theorem for manifolds with boundary}, to appear in Inter. Jour. of Geom. Meth. in Mode. Phys.\\
 \noindent [Wo] M. Wodzicki,  {\it Local invariants of
spectral
asymmetry}, Invent.Math. 75 no.1 143-178, 1984.\\

 \indent{  School of Mathematics and Statistics,
Northeast Normal University, Changchun Jilin, 130024 , China }\\
\indent E-mail: {\it wangy581@nenu.edu.cn; }\\

\end{document}